\documentclass[11pt,a4paper,reqno]{amsart}
\usepackage{graphicx}
\usepackage{epsfig}
\usepackage{color}
\usepackage{cite}
\usepackage{amssymb,euscript,appendix, lipsum}
\usepackage{listings,xcolor}
 \theoremstyle{definition}
 
 \theoremstyle{remark}
 \newtheorem{rem}{Remark}
\newcommand{\dist}{dist}
\begin{document}
\setcounter{page}{1}
\bigskip
\bigskip
\title[Nicolae Tarfulea: Functional Relationships in Solutions to ODE Systems]{Sparse Discovery of Functional Relationships in Solutions to Systems of Differential Equations}
\author[]{Nicolae Tarfulea$^{1,*}$}
\thanks{$^1$Department of Mathematics and Statistics, Purdue University Northwest, USA, \,e-mail: ntarfule@purdue.edu \\ \indent$^*$Corresponding author e-mail: ntarfule@purdue.edu}

\begin{abstract}
This work develops a framework to discover relations between the components of the solution to a given initial-value problem for a first-order system
of ordinary differential equations. This is done by using sparse identification techniques on the data represented by the numerical solution of the initial-value problem at hand.  
The only assumption is that there are only a few terms that connects the components, so that the mathematical relations to be discovered are sparse in the 
set of possible functions. We illustrate the method through examples of applications.   

\bigskip



\end{abstract}
\maketitle

\smallskip
\section{Introduction}
In this paper, we address the discovery of relationships between the solution components to a given system of ordinary differential equations (ODEs).  
Such interrelations in component solutions  might be elusive and difficult to discover, especially for phenomena modeled by large systems of 
nonlinear differential equations. Nonetheless, finding them could be very useful and involving computers in such
investigations leads to a renaissance in extracting patterns that are usually beyond human ability to grasp.  

This work is inspired by the quest to determine the underlying structure of nonlinear dynamical systems from data through the sparse identification of nonlinear dynamics (SINDy) methods. 
SINDy is a robust and versatile framework for uncovering the underlying dynamics of complex systems using sparse regression techniques, making it a valuable tool in scientific research and engineering applications. 
Developed in \cite{BPK}, SINDy leverages the principles of sparse regression to identify the simplest model $\dot{\mathbf{x}}= f(\mathbf{x})$, where $\mathbf{x}$ is the state vector and  
$f(\mathbf{x})$ represents the governing equations,  that describes the observed dynamics, making it particularly useful in scenarios where the underlying equations are unknown or only partially known. 
Although the aim of the work in this article is different, some ideas and techniques described here share some similarity to SINDy. Therefore, it is logical to begin by outlining the fundamental stages of SINDy.  
\begin{enumerate}
\item \textit{Data collection.} Collect time-series data of the system's states. This data could come from observations and/or experiment measurements. 
\item \textit{Library of Candidate Functions.} Construct a library, $\Theta (\mathbf{x})$, of candidate functions that might represent the system’s dynamics. 
This library might include polynomials, trigonometric functions, exponentials, or other basis functions that are hypothesized to describe the system. 
The method's success heavily depends on the choice of the candidate function library. 
\item \textit{Sparse Regression.} Use sparse regression techniques to identify which candidate functions from the library contribute significantly to the dynamics. 
These techniques enforce sparsity, identifying the most relevant terms in the library that govern the dynamics, as the goal is to find the model with the fewest terms  that accurately describes the data. 
Mathematically, the problem is to solve for the sparsest vector of coefficients, $\xi$, such that $ f(\mathbf{x})\approx \Theta (\mathbf{x})\xi$, where $\Theta (\mathbf{x})$ is the candidate library. 
The resulting sparse vector $\xi$ indicates which terms are significant, providing $f(\mathbf{x})$. 

\end{enumerate} 

Since its development, SINDy has seen numerous follow-ups and contributions in the literature. These contributions have expanded its application, improved its robustness, and integrated it with other methodologies. For example, Rudy, Brunton, Proctor, and Kutz \cite{RBPK} have adapted the SINDy framework to identify partial differential equations (PDEs) governing spatiotemporal data. Kang, Liao, and Liu \cite{KLL} have proposed techniques for identifying PDEs with numerical time evolution. 
They utilize Lasso for efficiency, a performance guarantee is established based on an incoherence property, and the main contribution is to validate and correct the results by time evolution error. 
Schaeffer, Tran, and Ward \cite{STW} have enhanced SINDy by incorporating group sparsity techniques, enabling the identification of dynamical systems with bifurcations. The method has been shown to effectively identify both the system dynamics and critical bifurcation points, providing insights into the system's behavior under parameter changes.  Loiseau and Brunton \cite{LB} have introduced a constrained version of SINDy that incorporates physical constraints into the sparse regression process. By enforcing constraints such as energy conservation or symmetries, the method improves the robustness and physical validity of the identified models. Boninsegna, N\"{u}ske, and  Clementi \cite{BNC} have extended SINDy to handle stochastic systems, using sparse learning to identify stochastic differential equations from data.
Champion, Lusch, Kutz, and Brunton \cite{CLKB} have extended SINDy to discover optimal coordinate transformations that simplify the underlying dynamics, integrating machine learning techniques to identify these transformations. 
Hoffmann, Nageshrao, and Haller \cite{HNH} have combined SINDy with cluster-based methods for the identification and control of nonlinear systems. 
Zhang and Schaeffer \cite{ZS} have provided theoretical analysis of the convergence of the SINDy algorithm. Forootani, Goyal, and Benner \cite{FGB}  have
proposed integrating neural networks with SINDy to improve robustness.  
Since its inception, SINDy has been applied in many areas of science and technology, including fluid dynamics, biological systems, neuroscience, physics, engineering, and machine learning; see \cite{AO,BC,FGB,LB,LNB,MBPK,SDN,ZKFB}, among many others.  

The objective of this work differs from SINDy in that the governing equations are known, and the goal is to identify functional relationships between the components of the solutions. 
Another key distinction is that the data is generated by numerically solving the system over a temporal grid, instead of using data from measurements. 
The method developed here systematically searches for relations between the solution components of ODE systems via sparse regression.  

The rest of this article is organized as follows: Section~\ref{sec2} introduces the technique behind identifying relationships between solution components. 
Section~\ref{sec3} then details the search algorithm for sparse identification of conserved quantities (i.e., solution constraints), with the corresponding  MATLAB code available in a public GitHub repository
\cite{T}.  Section~\ref{sec4} provides two illustrative examples. 
The paper ends with conclusions and an outlook presented in Section~\ref{conclusion}. 
\section{Identification of Conserved Quantities}
\label{sec2}
Let $D$ be an open subset of $\mathbb{R}\times\mathbb{R}^n$, $n\geq 1$, and let 
$f:D\to\mathbb{R}^n$ be a function $f(t,\mathbf{x})$ that is continuous in $t$ and Lipschitz continuous in $\mathbf{x}$. 
Consider the following initial-value problem associated to a system of first-order ordinary differential equations
\begin{eqnarray}
\label{ivp}
\dot{\mathbf{x}}(t)& = & f(t,\mathbf{x}(t)),\ t> t_0,\\ \label{ivpid}
\mathbf{x}(t_0)& = & \mathbf{x_0},
\end{eqnarray}
where the vector $\mathbf{x}(t)=(x_1(t),\, x_2(t),\,\ldots,\, x_n(t))^T\in \mathbb{R}^n$ denotes the state of the system at time $t$ and \eqref{ivpid} specifies the
initial state.  
Under the above conditions on the function $f$ and if $(t_0,\mathbf{x_0})\in D$, by Picard's local existence and uniqueness theorem, it is well-known that \eqref{ivp}-\eqref{ivpid} has a unique solution $\mathbf{x}(t)$ on a closed interval $[t_0, T]$, with $T>t_0$. 

In the present work, the goal is to find functional relationships  (or constraints) of the form 
\begin{equation}
\label{relation}
F(\mathbf{x}(t))=0,\quad t\in[t_0,T],
\end{equation}
where $\mathbf{x}(t)$ is the solution to \eqref{ivp}-\eqref{ivpid}. The form of the function $F$ will be specified in what follows. 

The discovery of such relations contributes to
the understanding of correlations among different components of the solution, and to find often hard to detect conserved quantities, supporting manifolds, and invariants. 
The search methodology is inspired by the SINDy methods and is presented next.

Following the approach introduced in \cite{BPK}, we first construct an augmented library $\Theta (\mathbf{x})$ consisting of candidate functions of the components 
of $\mathbf{x}$. There is truly a large freedom of choice in choosing the components of $\Theta (\mathbf{x})$, and 
the choice of candidate functions in the augmented library $\Theta (\mathbf{x})$ may not be
always clear. However, basic understanding of the mathematical model \eqref{ivp}-\eqref{ivpid} may provide good insight for 
a reasonable choice of the elements of $\Theta (\mathbf{x})$. 
It may contain a large number of elements, such as polynomial, exponential, logarithmic, and trigonometric terms, e.g., 
\begin{equation}
\label{theta}
\Theta (\mathbf{x})=\begin{bmatrix}
\mathbf{1} & \mathbf{x} & \mathbf{x^{(2)}} &\cdots &\mathbf{x^{(p)}}&\exp{(\mathbf{x})}&\ln{(\mathbf{x})}&\sin(\mathbf{x})&\cos(\mathbf{x})&\cdots
\end{bmatrix},
\end{equation}
where each of the entries denotes a specific row of functions in the components of $\mathbf{x}$. 
For example, here $\mathbf{x^{(2)}}=(x_1^2,\ x_1x_2, \  \ldots ,\ x_2^2,\ x_2x_3,\ \ldots,\ x_n^2)$ denotes the row vector of $2^{nd}$-order monomials 
formed with the components of $\mathbf{x}$ and
$\sin(\mathbf{x})=(\sin(x_1),\ \sin(x_2),\ \ldots,\ \sin(x_n))$.  

Each entry of $\Theta (\mathbf{x})$ represents a candidate component of $F$. 
That is, if $N$ is the length of the row vector-function $\Theta (\mathbf{x})$, we are looking for functions $F$ of the form:
$$F(\mathbf{x})=\sum_{i=1}^N\xi_i\Theta_i(\mathbf{x}),$$
where $\xi_i$ and $\Theta_i(\mathbf{x})$ are the constant coefficient and corresponding $i^{th}$-component of $\Theta (\mathbf{x})$, respectively, for each $i=1,\, 2,\,\ldots,\, N$. 
Henceforth, the coefficients $\xi_i$, $i=1,\, 2,\,\ldots,\, N$, will collectively be referred as ``$\xi$-coefficients.'' 

It is reasonable, although not necessary,  to assume that the function $F$ consists of only a few terms, 
making its composition sparse in the space generated by the augmented library $\Theta (\mathbf{x})$, that is, 
\begin{equation}
\label{F}
F(\mathbf{x})=\xi_{i_1}\Theta_{i_1}(\mathbf{x})+\xi_{i_2}\Theta_{i_2}(\mathbf{x})+\cdots+\xi_{i_k}\Theta_{i_k}(\mathbf{x}),
\end{equation}
with $1\leq i_1<i_2<\cdots<i_k\leq N$ and $k\ll N$.   
This leads to the search for $F$ as a sparse constant coefficient linear combination of the elements of $\Theta (\mathbf{x})$.  

That is, our goal is to set up a sparse-row regression problem to eliminate the zero (or close to zero) $\xi$-coefficients, which would lead to the finding of the nontrivial ones in the expression of $F$. 

\section{Algorithm}
\label{sec3}
As just mentioned, the search for the only few \textit{active} functions of the augmented library $\Theta (\mathbf{x})$ reduces to the finding of the nontrivial coefficients $\xi$ in the expression of $F$ in \eqref{F}. Inspired by the algorithm presented in \cite{BPK} for sparse identification of nonlinear dynamics (SINDy), we propose a least-squares algorithm for identification and  eliminating of all $\xi$-coefficients that are zero or close to zero in  the discrete maximum norm associated to the selected grid. 
Let us describe the proposed algorithm step by step.
\medskip

\noindent
\textbf{Step 1.} Find a highly-accurate numerical approximation of the true solution to \eqref{ivp}-\eqref{ivpid} at the grid points
$t_1$, $t_2$, ... , $t_m$, with $t_0<t_1<t_2<\ldots<t_m\leq T$, and arrange it into one large $(m+1)\times n$ matrix: 
 $$\mathbf{x}_{\text{grid}}=\begin{bmatrix}
\mathbf{x}^T_0\\
\mathbf{x}^T_1\\
\vdots\\
\mathbf{x}^T_m\\
\end{bmatrix}
=\begin{bmatrix}
x_{1;0} & x_{2;0} & \cdots &x_{n;0}\\
x_{1;1} & x_{2;1} & \cdots &x_{n;1}\\
\vdots & \vdots & \ddots & \vdots\\
x_{1;m} & x_{2;m} & \cdots &x_{n;m}\\
\end{bmatrix}.$$
Here, the first row, $\mathbf{x}^T_0=(x_{1;0},\, x_{2;0},\,\ldots, x_{n;0})$, is given by the initial data \eqref{ivpid}, and  $\mathbf{x}^T_j=(x_{1;j},\, x_{2;j},\,\ldots, x_{n;j})$ represents the numerical solution of \eqref{ivp}-\eqref{ivpid} at $t=t_j$, $j= 1,\, 2,\,\ldots,\, m$, i.e.,
$x_{i;j}$ approximates $x_i(t_j)$, for $i=1,\, 2,\,\ldots,\, n$ and $j= 1,\, 2,\,\ldots,\, m$.

There exists a multitude of highly accurate numerical methods that can be used, including the higher-order Taylor methods, Runge-Kutta methods, and multistep methods such as the explicit 
Adams-Bashforth methods and the implicit Adams-Moulton methods. In this paper, we employ the MATLAB built-in function \texttt{ode45}, which is a single-step solver based on an explicit Runge-Kutta (4,5) formula. This function provides high accuracy in approximating the solution of a system of ordinary differential equations, thanks to its fourth-order method and adaptive step size control (see \cite{ode45} for more details). 
\medskip

\noindent
\textbf{Step 2.} Construct the augmented matrix $\Theta_{\text{grid}}$ whose each row is obtained by using the numerical solution $\mathbf{x}_{\text{grid}}$ in the vector-function $\Theta (\mathbf{x})$ expression, that is: 
\begin{equation}
\label{thetam}
\Theta_{\text{grid}}=\begin{bmatrix}
\mathbf{1} & \mathbf{x_\text{grid}} & \mathbf{x_{\text{grid}}^{(2)}} &\cdots &\mathbf{x_{\text{grid}}^{(p)}}&\exp{(\mathbf{x_{\text{grid}}})}&\ln{(\mathbf{x_{\text{grid}}})}&\sin(\mathbf{x_{\text{grid}}})&\cos(\mathbf{x_{\text{grid}}})&\cdots
\end{bmatrix},
\end{equation}
where, for example, 
 $$\mathbf{x_{\text{grid}}^{(2)}} =\begin{bmatrix}
x_{1;0}^2 & x_{1;0}x_{2;0} & \cdots &x_{2;0}^2& x_{2;0}x_{3;0}&\cdots &x_{n;0}^2\\
x_{1;1}^2 & x_{1;1}x_{2;1} & \cdots &x_{2;1}^2&x_{2;1}x_{3;1}&\cdots &x_{n;1}^2\\
\vdots & \vdots & \ddots & \vdots &\vdots &\ddots &\vdots\\
x_{1;m}^2 & x_{1;m}x_{2;m} & \cdots &x_{2;m}^2&x_{2;m}x_{3;m}&\cdots &x_{n;m}^2\\
\end{bmatrix}
$$
and 
$$\sin(\mathbf{x_{\text{grid}}})=\begin{bmatrix}
\sin(x_{1;0}) & \sin(x_{2;0}) &\cdots &\sin(x_{n;0})\\
\sin(x_{1;1}) & \sin(x_{2;1}) &\cdots &\sin(x_{n;1})\\
\vdots & \vdots & \ddots & \vdots\\
\sin(x_{1;m}) & \sin(x_{2;m}) &\cdots &\sin(x_{n;m})
\end{bmatrix}.
$$
As a generic notation, we use $\Theta_{\text{grid}}[j_1,\, j_2,\ldots,\, j_r]$ to denote the $(m+1)\times r$ submatrix of $\Theta_{\text{grid}}$ with the $j_1$-, $j_2$-, $\ldots$, $j_r$ - th columns. 
\medskip

\noindent
\textbf{Step 3.} Reduce $\Theta_{\text{grid}}$ by repeating the following procedure several, say $p$, times:
\begin{enumerate}
\item[1)] Compute the reduced row echelon form of the square matrix $\Theta_{\text{rref}}:=\Theta_{\text{grid}}^T\Theta_{\text{grid}}$. Two mutually exclusive possibilities arise:
\begin{enumerate}
\item[(a)]If $\Theta_{\text{rref}}$ is the identity matrix, \textbf{stop}: 
the linear system \eqref{lsps} admits only the ``trivial'' solution, implicating a failed search because the columns of $\Theta_{\text{grid}}$ are linearly independent. In this case, a different augmented library $\Theta(\mathbf{x})$ should be considered. 
\item[(b)]If $\Theta_{\text{rref}}$ is not the identity matrix, record the indices of its row vectors with only one nonzero entry, which is also the leading entry (or pivot). Obviously, these indices label the zero components of $\xi$, which we discard.  
\end{enumerate}

\item[2)] Form the new matrix $\Theta_{\text{grid}}$ by discarding the columns corresponding to the indices found in 2) (b).
\end{enumerate}

Let us use the same notation, $\Theta_{\text{grid}}$, for the outcome, $\Theta_{\text{grid}}[i_1,\, i_2,\ldots,\, i_k]$, of this step. Obviously, the indices $ i_1<i_2<\cdots<i_k$ correspond to the 
survivor columns of  the original $\Theta_{\text{grid}}$. 
\medskip

\noindent
\textbf{Step 4.}  Solve the least-squares problem
\begin{equation}
\label{lsp}
\Theta_{\text{grid}}\xi =\eta, 
\end{equation}
where $\xi$ denotes the unknown column vector in $\mathbb{R}^k$ and $\eta$ represents the approximation error due to the numerical method used to find the numerical solution. For high-order methods, we treat $\eta$ as negligible,
that is, $\eta=0$.  A least-squares solution of \eqref{lsp} is a vector $\hat\xi$ in $\mathbb{R}^k$ such that 
$\dist (0, \Theta_{\text{grid}}\hat\xi)\leq \dist (0, \Theta_{\text{grid}}\xi)$ for all other vectors $\xi$ in $\mathbb{R}^k$. Here, $\dist $ denotes the Euclidian distance between vectors. 
Preserving the same notation, $\xi$, for the unknown, the set of least-squares solutions of \eqref{lsp} coincides with the nonempty set of solutions of the linear system
\begin{equation}
\label{lsps} 
\Theta_{\text{grid}}^T\Theta_{\text{grid}}\xi =0, 
\end{equation}
where $\Theta_{\text{grid}}^T$ denotes the transpose of the matrix $\Theta_{\text{grid}}$. 

The non-trivial solution set $\{ \xi_{i_1},\, \xi_{i_2},\, \ldots,\, \xi_{i_k}\}$ of the system \eqref{lsps} is what we are looking for. Here, the indices $\{ i_1,\, i_2,\, \ldots,\, i_k\}$ are associated to the 
corresponding candidate functions. 
\smallskip

A summary of the algorithm is given below; the MATLAB code in \cite{T} closely follows this algorithm.

\noindent\rule{\textwidth}{1.0pt}

\noindent\textbf{Algorithm}  Functional Relationships in Solutions to ODE Systems

\noindent\rule{\textwidth}{0.4pt}
 
\noindent\textbf{Input:} The initial-value problem (IVP) \eqref{ivp}-\eqref{ivpid}, the end-value $T$, the number of equally-spaced grid points $m$, and the number of iterations $p$.

\smallskip

\noindent\textbf{[Step 1]} Find a numerical approximation $\mathbf{x}_{\text{grid}}$ of the IVP solution at the grid points. 

\smallskip

\noindent\textbf{[Step 2]} Construct the augmented matrix $\Theta_{\text{grid}}$ as in \eqref{thetam}.

\smallskip

\noindent\textbf{[Step 3]} Refine $\Theta_{\text{grid}}$ $p$ times as follows: 
\begin{itemize}
\item Compute the reduced row echelon form of  $\Theta_{\text{rref}}:=\Theta_{\text{grid}}^T\Theta_{\text{grid}}$. 
\item If $\Theta_{\text{rref}}$ is the identity matrix, \textbf{stop} - a different augmented library should be considered. 
\item Otherwise, update $\Theta_{\text{grid}}$ by discarding the columns with just one nonzero entry and start over again. 
\end{itemize}
\smallskip

\noindent\textbf{Output:} The $\xi$-coefficients, which is the least-square solution to $\Theta_{\text{grid}}\xi=0$.  

\noindent\rule{\textwidth}{1.0pt}

\begin{rem}
As expected, the outcome of the algorithm may depend on the initial data  \eqref{ivpid}. Interestingly enough, one can then conclude that certain quantities involving the components of solutions to \eqref{ivp} are invariant, which could lead to simplifications and other interesting implications in the analysis of the system \eqref{ivp}. The applications in the next section illustrate this fact. 
\end{rem}
\section{Applications}
\label{sec4}
We apply the method outlined in Section~\ref{sec3} to two real-world scenarios. The first example demonstrates the method using a relatively straightforward system derived from a mathematical model of enzyme dynamics. 
In contrast, the second example involves a more intricate system: a mathematical model of oscillatory behavior in yeast glycolysis, incorporating the concentrations of seven biochemical species.

\subsection{Enzyme Dynamics} 
\label{ed}
\begin{figure}[h] % 'h' means place here
    \centering
    \includegraphics[width=1.1\textwidth]{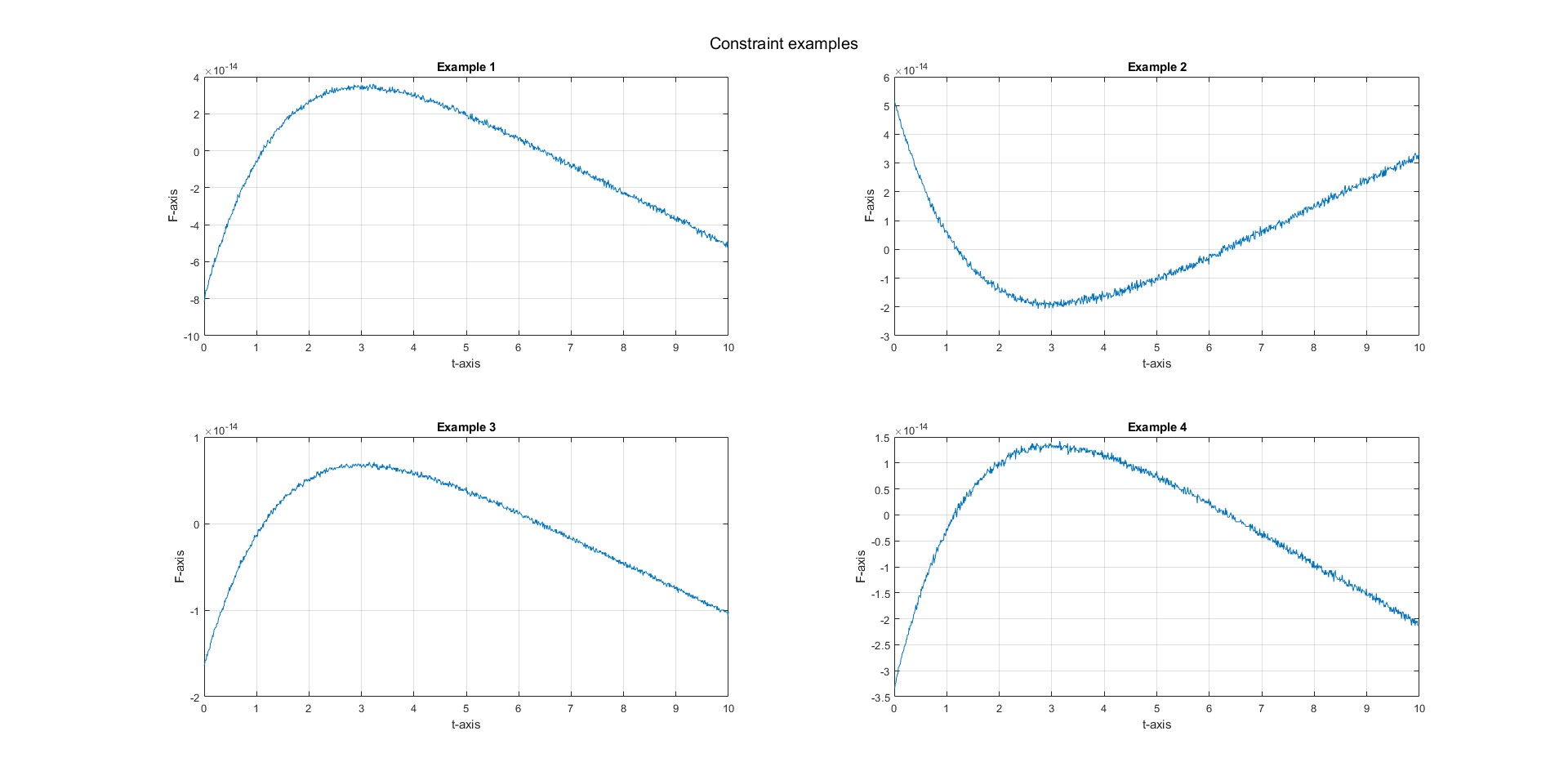} % adjust width
    \caption{\textbf{Example 1:} $\xi_1=12$, $\xi_2=-2$, $\xi_3=3$, $\xi_4=-5$, $\xi_5=-5$\newline
\textbf{Example 2:} $\xi_1=-14$, $\xi_2=8$, $\xi_3=5$, $\xi_4=3$, $\xi_5=3$\newline
\textbf{Example 3:} $\xi_1=3$, $\xi_2=-1$, $\xi_3=0$, $\xi_4=-1$, $\xi_5=-1$\newline 
\textbf{Example 4:} $\xi_1=7$, $\xi_2=-3$, $\xi_3=-1$, $\xi_4=-2$, $\xi_5=-2$}
    \label{1x}
\end{figure}

Let us consider one of the mathematical models for enzyme dynamics presented  in \cite[Chapter 2]{F}.  
Assume two species of proteins, $P$ and $L$, interact to form a complex $Q$ at a rate $k_1$, while $Q$ breaks down to its 
components $P$ and $L$ at a rate $k_{-1}$. Here, $P$, $L$, and $Q$ are concentrations, with unit $g/cm^3$, the
reaction rate $k_1$ is taken in unit $cm^3/g\cdot day$, and $k_{-1}$ is taken in unit of $1/day$.  The law of mass action says that
\begin{equation*}
\frac{dP}{dt}=-k_1PL+k_{-1}Q,\ \frac{dL}{dt}=-k_1PL+k_{-1}Q,\text{ and } \frac{dQ}{dt}=k_1PL-k_{-1}Q.
\end{equation*}
In the presence of an enzyme $E$ (a catalyst that promotes reaction rates), the conversion of proteins $S$ (called substrate) to protein $P$ (called product)  is sped up. Suppose $S$ and $E$ interact to form a complex $C_1$ at 
the rate $k_1$, while $C_1$ breaks down to the components $S$ and $E$ at the rate $k_{-1}$. Furthermore, suppose $C_1$ breaks down
to the components $E$ and $P$ at a rate $k_2$. 
By the law of mass action, we obtain the following first-order system of differential equations (see \cite[Chapter 2]{F} for more information)
\begin{align}
\label{E1}
\frac{dC_1}{dt}&=k_1SE-(k_{-1}+k_2)C_1,\\ \label{E2}
\frac{dE}{dt}&=-k_1SE+(k_{-1}+k_2)C_1,\\ \label{E3}
\frac{dS}{dt}&=-k_1SE+k_{-1}C_1,\\ \label{E4}
\frac{dP}{dt}&=k_2C_1. 
\end{align}
Consider the following initial data for the system \eqref{E1}-\eqref{E4}
\begin{equation}
\label{id}
 C_1(0)=1, E(0)=0, S(0)=1, P(0)=1.
\end{equation}
Running the MATLAB code in \cite{T}, based on the algorithm described in Section~\ref{sec3} with $\Theta (\mathbf{x})=[\mathbf{1} \  \mathbf{x}]$ in \eqref{theta}, reveals the following connection  between the components of the solution of the 
initial-value problem \eqref{E1}-\eqref{id}:
\begin{equation}
\label{con1}
\xi_1\cdot 1+\xi_2\cdot C_1(t)+\xi_3\cdot E(t)+\xi_4\cdot S(t)+\xi_5\cdot P(t)=0,\quad\text{for all } t\geq 0,
\end{equation}
for $\xi_1=-\xi_3-3\xi_5$, $\xi_2=\xi_3+\xi_5$, and $\xi_4=\xi_5$, with $\xi_3$ and $\xi_5$  as free parameters.  
Figure~\ref{1x} illustrates the preservation of contraint \eqref{con1} for values of $\xi$ specified in the caption.

\begin{figure}[h] % 'h' means place here
    \centering
    \includegraphics[width=1.1\textwidth]{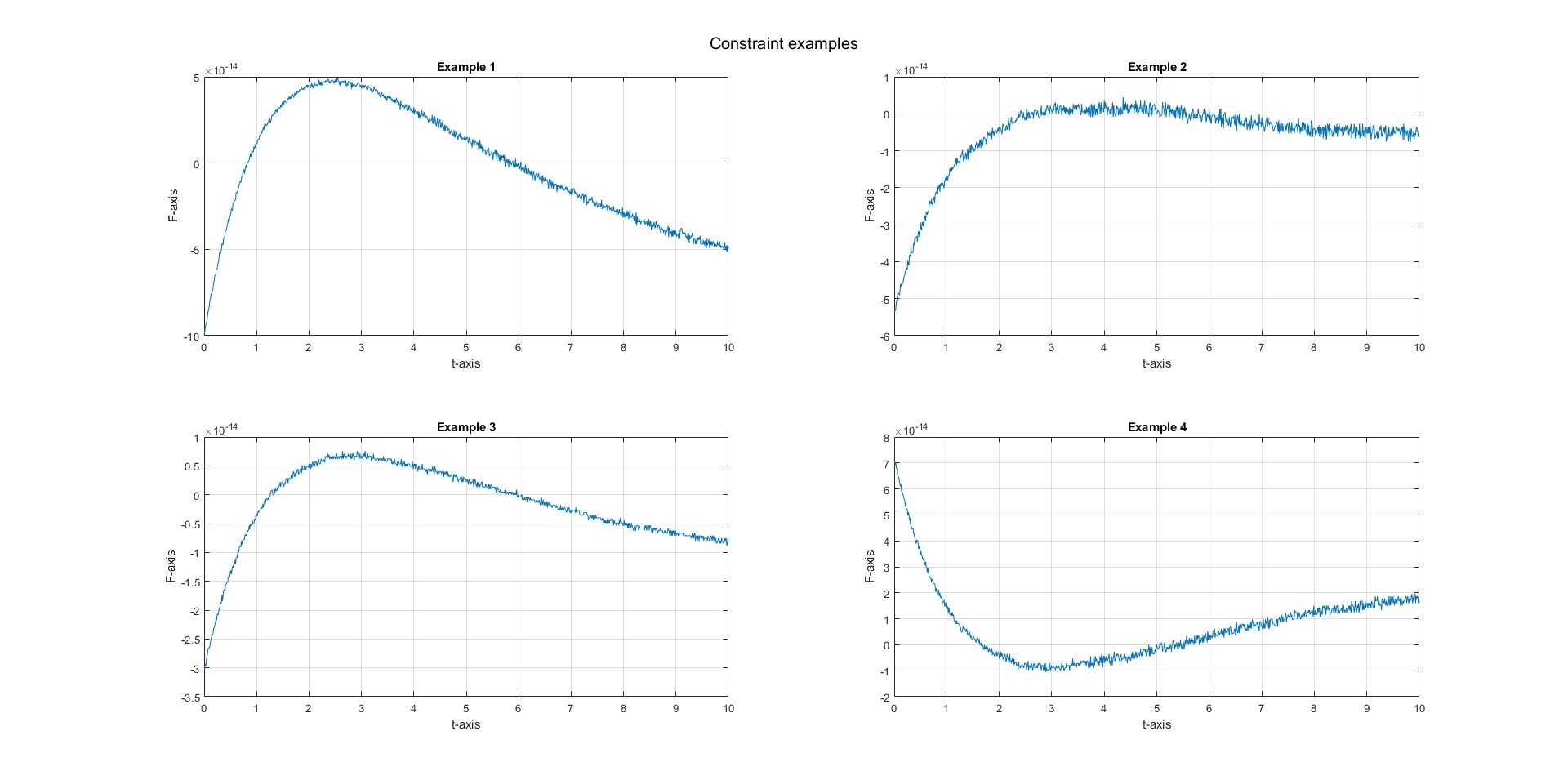} % adjust width
    \caption{\textbf{Example 1:} $\xi_1=-12$, $\xi_2=3$, $\xi_3=-2$, $\xi_4=5$, $\xi_5=5$\newline
\textbf{Example 2:} $\xi_1=-20$, $\xi_2=5$, $\xi_3=4$, $\xi_4=1$, $\xi_5=1$\newline
\textbf{Example 3:} $\xi_1=-8$, $\xi_2=2$, $\xi_3=1$, $\xi_4=1$, $\xi_5=1$\newline 
\textbf{Example 4:} $\xi_1=20$, $\xi_2=-5$, $\xi_3=-3$, $\xi_4=-2$, $\xi_5=-2$}
    \label{1x-2}
\end{figure}

If the initial data is changed to 
\begin{equation}
\label{id2}
 C_1(0)=2, E(0)=2, S(0)=1, P(0)=1,
\end{equation}
then the relation \eqref{con1} is valid for $\xi_1=-4\xi_3-4\xi_5$, $\xi_2=\xi_3+\xi_5$, and $\xi_4=\xi_5$, with $\xi_3$ and $\xi_5$ as free parameters.  Figure~\ref{1x-2} illustrates the preservation of constraint \eqref{con1} corresponding to the latter initial conditions for values of $\xi$ specified in the caption.

Notice that for the two sets of initial data considered above, only $\xi_1$ is different. It suggests that, for $\xi_2=\xi_3+\xi_5$, and $\xi_4=\xi_5$, with $\xi_3$ and $\xi_5$ as free parameters, 
$\xi_2\cdot C_1+\xi_3\cdot E+\xi_4\cdot S+\xi_5\cdot P$ is constant with respect to the independent variable $t$. This fact can easily be verified directly by noticing that
$$\frac{d}{dt}(\xi_2\cdot C_1+\xi_3\cdot E+\xi_4\cdot S+\xi_5\cdot P)=0.$$

\begin{figure}[h] % 'h' means place here
    \centering
    \includegraphics[width=1.1\textwidth]{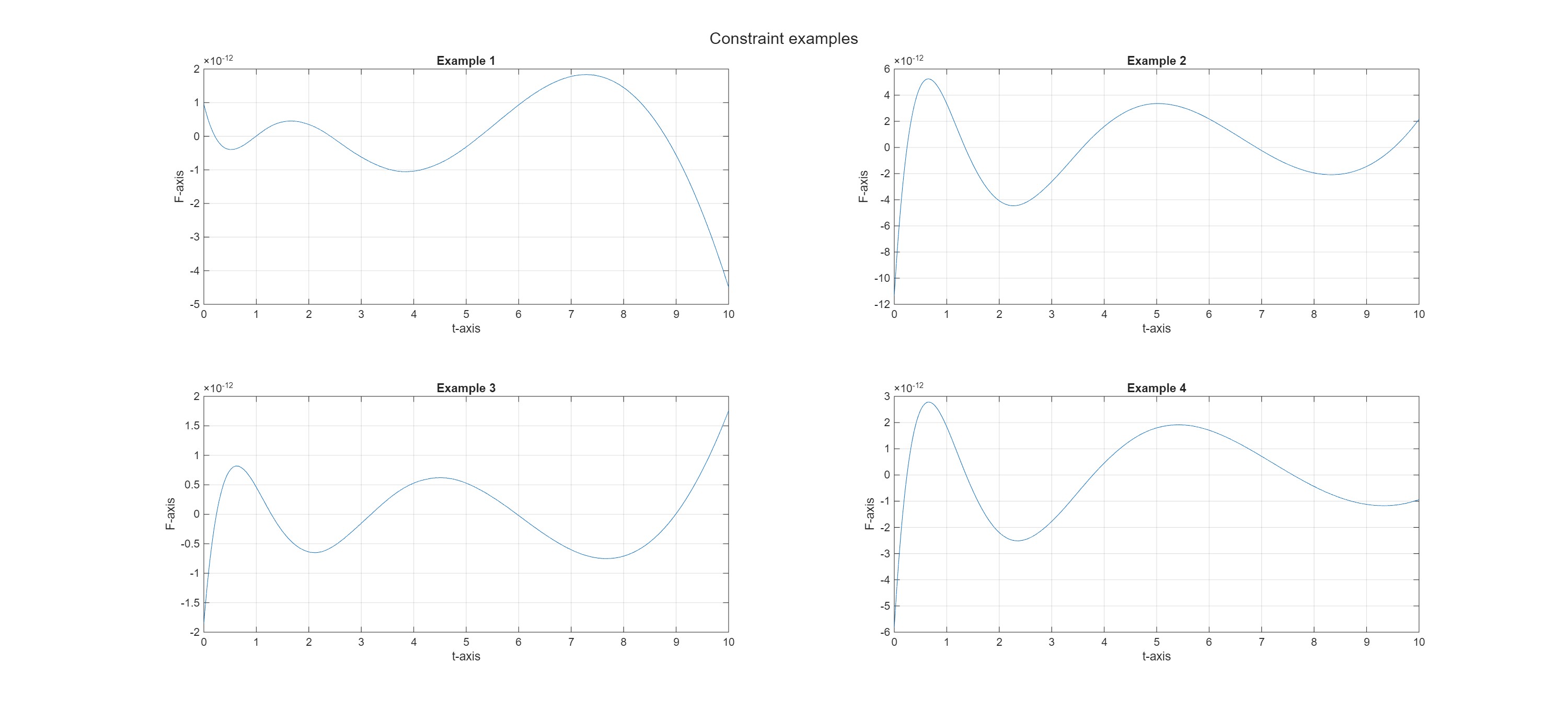} % adjust width
    \caption{\textbf{Example 1:} $\xi_1=27$, $\xi_2=-18$, $\xi_3=2$, $\xi_4=-25$, $\xi_5=-2$, $\xi_6=9$, $\xi_7=5$, $\xi_8=-5$, $\xi_9=4$, $\xi_{10}=-1$, $\xi_{11}=6$, $\xi_{12}=-1$, $\xi_{13}=3$, $\xi_{14}=3$, $\xi_{15}=-3$ \newline
\textbf{Example 2:}  $\xi_1=31$, $\xi_2=-9$, $\xi_3=0$, $\xi_4=-32$, $\xi_5=-1$, $\xi_6=0$, $\xi_7=2$, $\xi_8=2$, $\xi_9=9$, $\xi_{10}=3$, $\xi_{11}=6$, $\xi_{12}=-2$, $\xi_{13}=2$, $\xi_{14}=2$, $\xi_{15}=-4$\newline
\textbf{Example 3:} $\xi_1=-3$, $\xi_2=8$, $\xi_3=-4$, $\xi_4=11$, $\xi_5=0$, $\xi_6=1$, $\xi_7=5$, $\xi_8=-2$, $\xi_9=-8$, $\xi_{10}=1$, $\xi_{11}=-3$, $\xi_{12}=-3$, $\xi_{13}=3$, $\xi_{14}=-3$, $\xi_{15}=0$\newline 
\textbf{Example 4:} $\xi_1=21$, $\xi_2=-15$, $\xi_3=2$, $\xi_4=-25$, $\xi_5=4$, $\xi_6=6$, $\xi_7=5$, $\xi_8=1$, $\xi_9=5$, $\xi_{10}=-4$, $\xi_{11}=7$, $\xi_{12}=-4$, $\xi_{13}=-3$, $\xi_{14}=4$, $\xi_{15}=-3$ }
    \label{1x-3}
\end{figure}

For $\Theta (\mathbf{x})=[\mathbf{x}]$ and initial data \eqref{id}, the result is 
$$\xi_1\cdot C_1+\xi_2\cdot E+\xi_3\cdot S+\xi_4\cdot P=0,$$
for $\xi_1=-2\xi_4$, $\xi_2=-3\xi_4$, and $\xi_3=\xi_4$, with $\xi_4$ free parameter. Notice that the latter choice of $\Theta (\mathbf{x})$ uncovers only one invariant quantity
$$-2C_1-3E+S+P=0,$$
while the former finds a two parameter family of such invariant quantities. Unsurprisingly, a richer  augmented library $\Theta (\mathbf{x})$ leads to more numerous and interesting solution connections. 
The MATLAB code in \cite{T} produces interesting results for richer augmented libraries. For example, if $\Theta (\mathbf{x})=[\mathbf{1} \  \mathbf{x}\ \mathbf{x}^2]$,
then the result is 
\begin{align*}
\xi_1+\xi_2C_1+\xi_3E+\xi_4S+\xi_5P+\xi_6C_1^2+\xi_7C_1E+\xi_8E^2+\xi_9C_1S+\xi_{10}ES\\
+\xi_{11}S^2+\xi_{12}C_1P+\xi_{13}EP+\xi_{14}SP+\xi_{15}P^2=0,
\end{align*}
where $\xi_1=-\xi_3-3\xi_5-\xi_8-3\xi_{13}-9\xi_{15}$, $\xi_2=\xi_3+\xi_5-\xi_7+2\xi_8-3\xi_{12}+4\xi_{13}+6\xi_{15}$, $\xi_4=\xi_5-\xi_{10}+\xi_{13}-3\xi_{14}+6\xi_{15}$, 
$\xi_6= \xi_7-\xi_8+\xi_{12}-\xi_{13}-3\xi_{14}-\xi_{15}$, $\xi_9=\xi_{10}+\xi_{12}-\xi_{13}+\xi_{14}-2\xi_{15}$, and $\xi_{11}=\xi_{14}-\xi_{15}$, with
$\xi_3$, $\xi_5$, $\xi_7$, $\xi_8$, $\xi_{10}$, $\xi_{12}$, $\xi_{13}$, $\xi_{14}$, and $\xi_{15}$ free variables.
Figure~\ref{1x-3} illustrates the preservation of constraint \eqref{con1} corresponding to the initial conditions \eqref{id2} for values of $\xi$ specified in the caption.

\subsection{Glycolytic oscillator model} 
\label{gl}
We consider the model of oscillations in yeast glycolysis introduced in \cite{DN} (see also \cite{BPK,RCWH}). The model details are not critical to our purpose, we 
instead take this biological model as another example from which we want to extract information by using the algorithm described in Section~\ref{sec3}. The model consists of an ODE system for the concentrations of seven biochemical species:
\begin{align}
\label{gl1}
\frac{dS_1}{dt}&=J_0-\frac{k_1S_1S_6}{1+(S_6/K_1)^q},\\
\frac{dS_2}{dt}&=2\frac{k_1S_1S_6}{1+(S_6/K_1)^q}-k_2S_2(N-S_5)-k_6S_2S_5,\\
\frac{dS_3}{dt}&=k_2S_2(N-S_5)-k_3S_3(A-S_6),\\
\frac{dS_4}{dt}&=k_3S_3(A-S_6)-k_4S_4S_5-\kappa (S_4-S_7),\\
\frac{dS_5}{dt}&= k_2S_2(N-S_5)-k_4S_4S_5-k_6S_2S_5,\\
\frac{dS_6}{dt}&=-2\frac{k_1S_1S_6}{1+(S_6/K_1)^q}+2k_3S_3(A-S_6)-k_5S_6,\\ \label{gl7}
\frac{dS_7}{dt}&=\psi\kappa (S_4-S_7)-kS_7,
\end{align}
where the model parameters (with discarded units) are taken from \cite{DN}: $J_0=2.5$, $k_1=100$, $k_2=6$, $k_3=16$,
$k_4=100$, $k_5=1.28$, $k_6=12$, $k=1.8$, $\kappa=13$, $q=4$, $K_1=0.52$, $\psi=0.1$, $N=1$, and $A=4$.
The initial conditions for $S_1,\,\ldots,\  S_7$ will be chosen from the ranges provided in \cite{DN}, that is, 
$[0.15,\, 1.60]$,  $[0.19,\, 2.16]$, $[0.04,\, 0.20]$, $[0.10,\, 0.35]$, $[0.08,\, 0.30]$, $[0.14,\, 2.67]$, and $[0.05,\, 0.10]$, respectively. 

We tried different libraries of candidate functions $\Theta$ mentioned in \eqref{theta}. 
Running the MATLAB code in \cite{T} with $\Theta (\mathbf{x})=[\mathbf{1} \  \mathbf{x}]$, $\Theta (\mathbf{x})=[\mathbf{1} \  \mathbf{\sin (x)}]$, and  
$\Theta (\mathbf{x})=[\mathbf{1} \ \mathbf{x}\ \mathbf{\sin (x)}]$ shows that there
are no nontrivial functional relationships between the components of the solutions.   
For $\Theta (\mathbf{x})=[\mathbf{1} \  \mathbf{x}\ \mathbf{x^2}]$, the output consists of nontrivial linear combinations of 36 monomials 
of degrees zero, one, and two: 

\noindent
$\xi_1 +\xi_2 S_1+\xi_3 S_2+\xi_4 S_3+\xi_5 S_4+\xi_6 S_5+\xi_7 S_6+\xi_8 S_7+\xi_9 S_1^2+\xi_{10}S_1S_2+\xi_{11}S_2^2+\xi_{12}S_1S_3 
+\xi_{13}S_2S_3+\xi_{14}S_3^2+\xi_{15}S_1S_4+\xi_{16}S_2S_4+\xi_{17}S_3S_4+\xi_{18}S_4^2+\xi_{19}S_1S_5+\xi_{20}S_2S_5+\xi_{21}S_3S_5+\xi_{22}S_4S_5+\xi_{23}S_5^2+\xi_{24}S_1S_6
+\xi_{25}S_2S_6+\xi_{26}S_3S_6+\xi_{27}S_4S_6+\xi_{28}S_5S_6+\xi_{29}S_6^2+\xi_{30}S_1S_7+\xi_{31}S_2^1S_7+\xi_{32}S_3S_7+\xi_{33} S_4S_7+\xi_{34}S_5S_7+\xi_{35}S_6S_7+\xi_{36}S_7^2 
= 0$. 

\begin{figure}[h] % 'h' means place here
    \centering
    \includegraphics[width=0.85\textwidth]{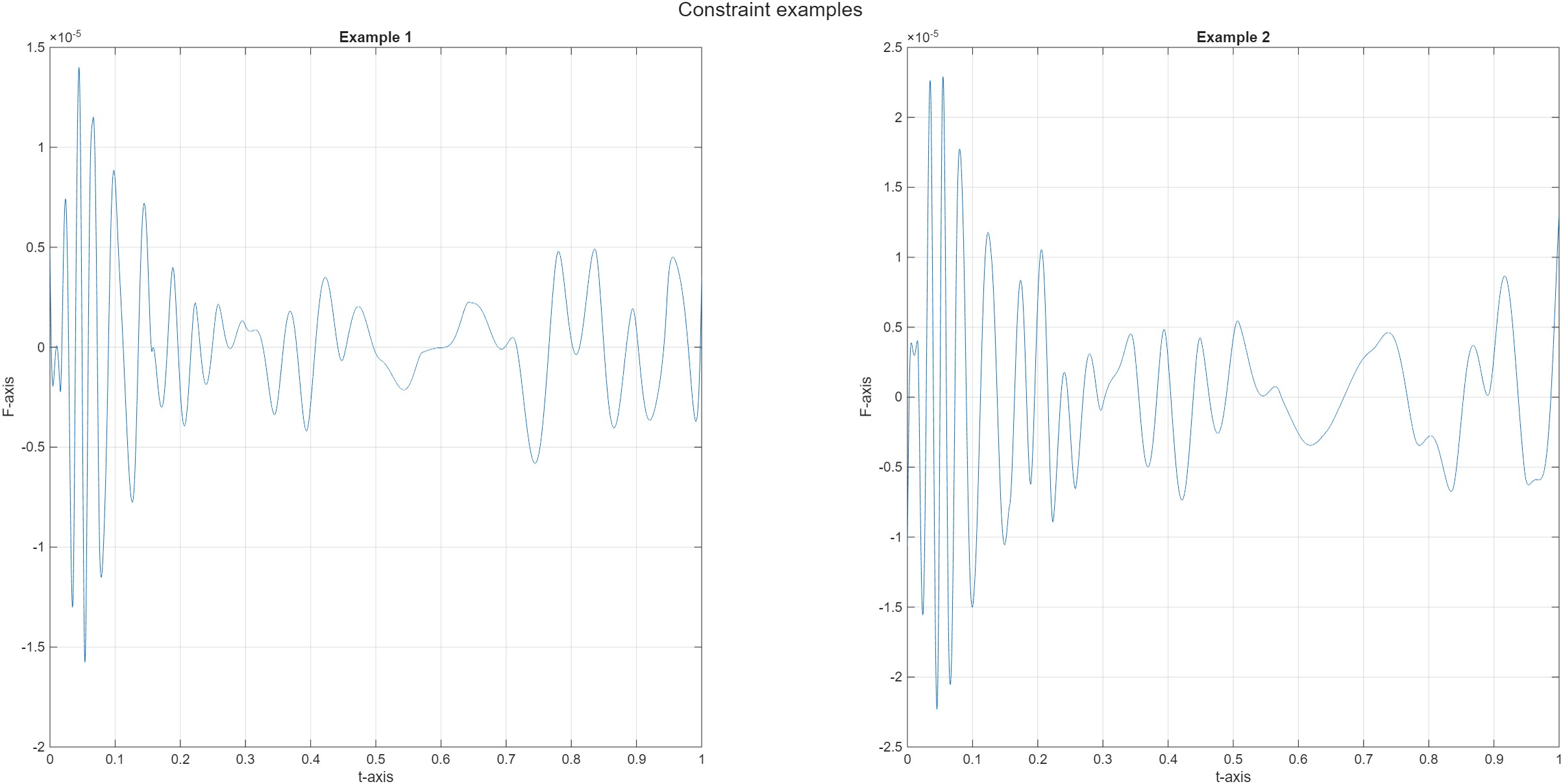} % adjust width
    \caption{\textbf{Example 1:} $\xi_1=3.7753$,   $\xi_2=-4.5229$,  $\xi_3= -9.2616$,  $\xi_4=15.6269$,  $\xi_5=17.2988$,  $\xi_6= 10.0912$,  $\xi_7=-1.0314$,  $\xi_8=6.5885$,  $\xi_9=1.0066$,  $\xi_{10}= 3.9460$,  $\xi_{11}=2.9008 $,  
$\xi_{12}=  -5.3707$,  $\xi_{13}=-0.4891$,  $\xi_{14}=-19.2035$,  $\xi_{15}=-5.5577$,  $\xi_{16}= -6.1561$,  $\xi_{17}= -26.0902$,  $\xi_{18}=-1.3775$,  $\xi_{19}=-2.8670$,  $\xi_{20}=-6.4380$,  $\xi_{21}= 10.0405$,  
$\xi_{22}=10.9240$,  $\xi_{23}= -4$,  $\xi_{24}=1.0466$,  $\xi_{25}=1.6645$,  $\xi_{26}=-3$,  $\xi_{27}=-4$,  $\xi_{28}=-3$,  $\xi_{29}=-0.0623$,  $\xi_{30}=-5$,  $\xi_{31}=1$,  $\xi_{32}=-2$,  $\xi_{33}=0$,  $\xi_{34}=2 $,  $\xi_{35}=0$,  $\xi_{36}=0$
\newline
\textbf{Example 2:} $\xi_1=-8.1704$,   $\xi_2=7.5115$,   $\xi_3= 8.6726$,   $\xi_4=16.7390$,   $\xi_5= 0.5772$,   $\xi_6=-17.0689$,   $\xi_7=1.5089$,   $\xi_8=-2.1254$,   $\xi_9= -1.5832$,   $\xi_{10}=-3.9060$,   $\xi_{11}=-1.9235$,   $\xi_{12}=-5.3588
$,   $\xi_{13}= -11.5336$,   $\xi_{14}=   -0.1188$,   $\xi_{15}=   -0.0902$,   $\xi_{16}=-0.1688$,   $\xi_{17}=14.3894$,   $\xi_{18}=-1.5993$,   $\xi_{19}=5.6994$,   $\xi_{20}=8.4104$,   $\xi_{21}=14.4073$,   $\xi_{22}=-9.5428$,   $\xi_{23}=-1$,   $\xi_{24}=-1.2154$,   
$\xi_{25}=-0.8209$,  $\xi_{26}=-4$,   $\xi_{27}=0$,   $\xi_{28}=4$,   $\xi_{29}=0.1935$,   $\xi_{30}=4$,   $\xi_{31}=-3$,   $\xi_{32}=-3$,   $\xi_{33}=1$,   $\xi_{34}=2$,   $\xi_{35}=-1$,   $\xi_{36}= -3$}
\label{2x-4}
\end{figure}

As an example, for the system \eqref{gl1}-\eqref{gl7} endowed with the initial data 
\begin{equation}
\label{idgl}
S_1(0)=1,\ S_2(0)=1, \ S_3(0)=0.1,\ S_4(0)= 0.2,\ S_5(0)= 0.2,\ S_6(0)= 1,\ S_7(0)= 0.1, 
\end{equation}
the algorithm finds 
a large solution involving 25 basic and 11 free $\xi$-coefficients. Since the output is too long to be displayed here, the interested reader is invited to run the code in \cite{T} to see the full solution. Figure~\ref{2x-4} 
 illustrates the preservation of the above constraint corresponding to the initial data \eqref{idgl} for values of $\xi$ specified in the caption.

Other nontrivial results can be obtained for various choices of libraries of candidate functions, e.g., $\Theta (\mathbf{x})=[\mathbf{1} \  \mathbf{x} \ \mathbf{\sin (x)}\ \mathbf{\sin(2x)}]$.
It is important to remember that the results are in the least-squares sense, meaning that the $\xi$-coefficients are the best one can get for a specific library of candidate functions $\Theta$.

\section{Conclusion} 
\label{conclusion}
In this work, we propose and demonstrate a novel technique for identifying functional relationships within the solutions of systems of ordinary differential equations (ODEs). 
Such conserved quantities are often challenging to uncover, particularly in applications involving large nonlinear systems. 
Identifying these invariants can significantly enhance both qualitative and quantitative analyses, potentially revealing insightful patterns in the solution behavior. 
Leveraging computational tools in this context greatly improves the efficiency and effectiveness of detecting such patterns, which might otherwise remain hidden.

The proposed method is inspired by the well-established Sparse Identification of Nonlinear Dynamics (SINDy) framework. 
The core idea involves selecting the most relevant terms from a predefined library of candidate functions using sparse, nonlinear regression techniques. 
In our approach, the data (namely, numerical solutions of systems of ODEs) are generated using standard numerical solvers.

We illustrate the method using two examples: an enzyme kinetics model and a glycolytic oscillator model. A key contribution of this work is a MATLAB function \cite{T} that implements the proposed algorithm. 
The code is designed with sufficient generality to be applicable to arbitrary systems of ODEs.

Several promising directions for future research emerge from this study. One avenue involves extending the method to accommodate functional relationships with non-constant $\xi$-coefficients. 
Another compelling direction is the generalization of the approach to systems of partial differential equations (PDEs), which introduces additional layers of complexity, including multidimensional data 
generation and the construction of higher-dimensional spaces of candidate functions. 

\bigskip

\noindent
\textbf{Acknowledgments} This work was partially supported by a Purdue University Northwest Catalyst Grant.

\end{document}